\documentclass[a4paper,11pt]{article}

\title{The intrinsic fundamental group of a linear category}

\author{Claude Cibils, Mar\'\i a Julia Redondo and Andrea Solotar
\thanks{This
work has been supported by the projects MATHAMSUD, UBACYTX212, PIP-CONICET
112-2000801-00487 and PICT-2007-02182 (ANPCyT). The second and third authors are research members of CONICET (Argentina) and the
third author is a Regular Associate of ICTP Associate Scheme.}
}

\date{}

\input{xy}
\xyoption{all}
\xyoption{poly}
\usepackage[all]{xy}
\usepackage{amsmath,amsthm,amsfonts,amssymb,stmaryrd}
\usepackage{latexsym}
\usepackage{color,palatino}
\usepackage{graphicx}%poner [dvips] para gr\'{a}ficos
\usepackage{float}  %Este packete permite controlar la flotacion de figuras al poner [H] o [h] como opcion de \begin{figure}
\RequirePackage{amssymb}
\RequirePackage[T1]{fontenc}

\usepackage[
  breaklinks,
  naturalnames,
  ]{hyperref}

\newtheorem{thm}{Theorem}[section]
\newtheorem{cor}[thm]{Corollary}
\newtheorem{lem}[thm]{Lemma}
\newtheorem{pro}[thm]{Proposition}
\newtheorem{defi}[thm]{Definition}%[section]
\newtheorem{rem}[thm]{Remark}%[section]
\newtheorem{exa}[thm]{Example}

\let\oldqed\qed
\renewcommand\qed{\oldqed\par\bigskip}

\def\A{{\mathcal A}}
\def\B{{\mathcal B}}
\def\C{{\mathcal C}}
\def\D{{\mathcal D}}

\def\K{{\mathcal K}}

\def\U{{\mathcal U}}

\def\ker{{\mathrm {ker}}}

\def\Hom{{\mathrm {Hom}}}
\def\aut{{\mathrm {Aut}}}
\def\Aut{{\mathrm {\bf Aut}}}

\def\charac{{\mathrm {char}}}

\def\lim{\mathop{\rm lim}\nolimits}

\def\st{\mathsf{St}}

\def\gal{\mathsf{Gal}}
\def\cov{\mathsf{Cov}}

\begin{document}

\sf

\maketitle

\begin{abstract}
We provide an intrinsic definition of the fundamental group of a li\-near category over a
ring as the automorphism group of the fibre functor on Galois coverings. If the universal covering exists, we prove that
this group is isomorphic to the Galois group of the universal covering.
The grading deduced from a Galois covering enables us to describe
the canonical monomorphism from its automorphism group to the first Hochschild-Mitchell
cohomology vector space.
\end{abstract}


\section{\sf Introduction}

The purpose of this work is to provide a positive answer to the question of the existence
of an intrinsic and canonical fundamental group $\pi_1$ associated to a $k$-category
$\B$, where $k$ is a commutative ring. {The fundamental group we introduce takes into
account the linear structure of the category $\B$, it differs from the fundamental group
of the underlying category obtained as the classifying space of its nerve
(\cite{se,qu,sc}).}

The fundamental group that we define is intrinsic in the sense that it does not depend on
the presentation of the $k$-category by generators and relations. {In case a universal
covering exists, we obtain that the fundamental groups {constructed by R.
Mart\'inez-Villa and J.A. de la Pe\~{n}a (see \cite{MP}, and \cite{boga,ga,asde})} depending
on a presentation of the category by a quiver and relations are in fact quotients of the
intrinsic $\pi_1$ that we introduce. Note that those groups can vary according to
different presentations of the same $k$-category (see for instance \cite{asde,buca,le1})
while the group that we introduce is intrinsic, since we define it as the automorphisms
of the fibre functor of the Galois coverings over a fixed object.}

In fact if a universal covering $U: \U\rightarrow\B$ exists, the fundamental group that
we define is isomorphic to the automorphism group $\mathsf{Aut}U$, and in this case
changing the base object provides isomorphic intrinsic fundamental groups.

The methods we use are inspired in the topological case considered for instance in R.
Douady and A. Douady's book \cite{dodo}. They are closely related to the way in which the
fundamental group is considered in algebraic geometry after A. Grothendieck and C.
Chevalley.

In algebraic topology a space has a universal cover if it is connected, locally path-connected and semi-locally simply connected.
In other words, usually a space has a universal cover. By contrast, linear categories do not have
in general universal coverings.

Our work is very much indebted to the pioneer work of P. Le Meur in
his thesis \cite{le2}, see also \cite{le}. He has shown that under some hypotheses on the
category, there exists an optimal fundamental group in the sense that all other
"fundamental groups" deduced from different presentations are quotients of the optimal
one. {His method consists mainly in tracing all the possible presentations of a given
category, and relating the diverse resulting "fundamental groups". As already quoted, we
adopt a different point of view in this paper.

In Section \ref{coverings} we recall the definition of a covering of a $k$-category and
we prove properties about morphisms between coverings as initiated in \cite{le,le2}. In
Section \ref{galoiscoverings} we define Galois coverings and next we study some
properties of this kind of coverings.  The main results are Theorem \ref{structure}, which
describes the structure of Galois coverings and Theorem \ref{morphism2}, which concerns
morphisms between Galois coverings and the relation between the associated groups of
automorphisms. We provide the definition of the universal covering in the category of
Galois coverings of a fixed $k$-category $\B$. In a forthcoming paper we will study the
behaviour of Galois coverings through fibre products, as well as a criterion for a
covering to be Galois or universal. Differences with the usual algebraic topology setting
will also appear, since the fibre product of coverings of $k$-categories does not provide in
general a covering through the projection functor.

In Section \ref{fundamental} we define
the intrinsic fundamental group $\pi_1(\B, b_0)$ and we prove some properties of this new object. If the universal covering exists, we prove that this group is isomorphic to the Galois group of the universal covering. In \cite{crs} we provide explicit computations of the intrinsic fundamental group of some algebras. In particular
we compute the fundamental group of $M_p(k)$, where $p$ is prime and $k$ an algebraically closed field of characteristic zero, which is the direct product of the free group on $p-1$ generators with the cyclic group of order $p$. The fundamental group of triangular matrix algebras is the free group on $n-1$ generators.
The fundamental group of the truncated
polynomial algebra $k[x]/(x^p)$ in characteristic $p$ is the product of the infinite cyclic group and the cyclic group of order $p$. In case $k$ is a field containing all roots of unity of order $2$ and $3$,
we prove that $\pi_1(k^3)=C_2\times C_3$, while if $k$ contains all roots of unity of order  $3$ and $4$,
we obtain that $\pi_1(k^4)= (C_2 * C_2)\times C_6 \times C_4\times C_2.$

In section \ref{fundamental} we also show that if the universal covering exists the fibre functor induces an equivalence between the category of Galois coverings of $\B$ and the subcategory of $\pi_1(\B, b_0)$-$\mathsf{Sets}$ whose objects are sets with a transitive action of the group $\pi_1(\B, b_0)$ such that the isotropy group of an element is invariant.

In the last section we suppose  that $k$ is a field and that the endomorphism algebra of each object
of the $k$-category  is reduced to $k$. We recover in a simple way the canonical
$k$-linear embedding (see \cite{asde,jap-s,cire}) of the abelian characters of the
automorphism group of a Galois covering to the first Hochschild-Mitchell cohomology
vector space of the category. As an immediate consequence, {if there exists a Galois covering whose group is isomorphic to the fundamental group,
the abelian characters of the fundamental group embed into the cohomology of degree one.}
For this, we use our description of Galois coverings as well as the
canonical  grading  of the $k$-category deduced from a Galois covering, as obtained in
\cite{cima}. In this way Euler derivations are considered, see also \cite{fagegrma,fagrma}.

Note that it will be interesting to explore the behaviour of the intrinsic fundamental
group with respect at least to Morita equivalences of $k$-categories. As expected the
fundamental group is an invariant of the equivalence class of a $k$-category but not of
its Morita class. Of course in case the category admits a unique basic representative in
its Morita class, the fundamental group attached to this category can be considered as
the canonical fundamental group of the Morita equivalence class.

{\textbf{Acknowledgements:}  We thank Alain Brugui\`{e}res, Jos\'{e} Antonio de la Pe\~na, Patrick Le Meur, Christian Pauly and Sonia Trepode
for useful discussions. {In particular Patrick Le Meur pointed out the importance of taking into account the automorphisms of the base category.} The first author wish to thank Prof. Pu Zhang and Jiao Tong
University for an invitation to Shanghai where part of this work has been performed.}

%%%%%%%%%%%%%%%%%%%%%%%%%%%%%%%%%%%%%%%%%%%%%%%%%%%%%%%%%%%%%%%%%%%%%

\section{\sf $k$-categories, stars and coverings}\label{coverings}
Let $k$ be a commutative ring. A \emph{$k$-category} is a small category $\B$ such that
each morphism set ${}_y\B_x$ from an object $x\in\B_0$ to an object $y\in\B_0$ is a
$k$-module, the composition of morphisms is $k$-bilinear and {$k$-multiples of the
identity at each object are central in the endomorphism ring of the object}. Note that
such $k$-categories are also called \emph{linear categories over $k$}. In particular each
endomorphism set of an object is a $k$-algebra, and ${}_y\B_x$ is a ${}_y\B_y -
{}_x\B_x$-bimodule.

Each $k$-algebra $A$ provides a single object $k$-category $\B_A$ with endomorphism ring
$A$. The structure of $A$ can be described more precisely by choosing a finite set $E$ of
orthogonal idempotents of $A$, such that $\sum_{e\in A}e=1$ in the following way: the
$k$-category $\B_{A,E}$ has set of objects $E$ and morphisms from $e$ to $f$ the
$k$-module $fAe$. Note that $\B_{A,\{1\}}=\B_A$. {This approach is meaningful since
clearly the category} of left $A$-modules is isomorphic to the category of $k$-functors
from $\B_{A,E}$ to the category of $k$-modules, where a $k$-functor is a functor which is
$k$-linear when restricted to morphisms.

\begin{defi} The star $\st_{b_0}\B$ of a $k$-category $\B$ at an
object $b_0$ is the direct sum of all the morphisms with source or target $b_0$ :
\[\st_{b_0}\B = \left(\bigoplus_{y\in\B_0} \ {}_y\B{_{b_0}}\right)\ \oplus
\ \left(\bigoplus_{y\in\B_0} \ {}_{b_0}\B{_y}\right).\] Note that this $k$-module counts
twice the endomorphism algebra at $b_0$.
\end{defi}

\begin{defi}\label{covering}
 Let $\C$ and $\B$ be $k$-categories. A $k$-functor
$F:\C\rightarrow\B$ is a covering of $\B$ if it is surjective on objects and if $F$
induces $k$-isomorphisms between the corresponding stars. More precisely, for each
$b_o\in\B_0$ and each $x$ in the non-empty fibre $F^{-1}(b_0)$, the map
\[F_{b_0}^x:\st_x\C\longrightarrow\st_{b_0}\B.\]
induced by $F$ is a  $k$-isomorphism.

\end{defi}

\begin{rem} Each star is the direct sum of the source star $\st^-_{b_0}\B =
\bigoplus_{y\in\B_0}{}_y\B{_{b_0}}$ and the target star
$\st^+_{b_0}\B=\bigoplus_{y\in\B_0} \ {}_{b_0}\B{_y}$. Since   $\st^-$ and $\st^+$ are
preserved under any $k$-functor, the condition of the definition is equivalent to the
requirement that the corresponding target and source stars are
isomorphic through $F$.\\
Moreover this splitting goes further :
for $b_1\in\B_0$, the restriction of $F$ to
$\bigoplus_{y\in F^{-1}(b_1)} {}_y\C_x$, for all $x \in F^{-1}(b_0)$
is $k$-isomorphic to the corresponding $k$-module
${}_{b_1}\B_{b_0}$. The same holds with respect to the target star and morphisms starting
at all objects in a single fibre.
\end{rem}

\begin{rem}
The previous facts show that Definition \ref{covering} coincides with the one given by K.
Bongartz and P. Gabriel in \cite{boga}.
\end{rem}

To each small category $\A$ one may associate its linearization $k\A$ in the following
way: objects of $k\A$ are the objects of $\A$, while morphisms are free $k$-modules on
the sets of morphisms of $\A$. Such linearized $k$-categories admit by construction a multiplicative basis of morphisms, which is not usually the case, see for instance \cite{bagarosa}.
Hence a $k$-category is not in general the linearization of a small category.
Note that any usual covering $F:\C \to \B$ of categories provides by linearization a covering $kF:k\C \to k\B$.

\begin{exa}\cite{le}\label{nongalois}
Consider the following $k$-categories $\C$ and $\K$ obtained by linearization of the
categories given by the corresponding diagrams ($\K$ is called the Kronecker category):

\[
\xymatrix@R-12pt {& & {s_0}
\ar[dl]_{\alpha_0}
\ar[dr]^{\beta_0}   &
\\
{\C} : &t_0 & & t_1
\\
&& s_1
\ar[ul]^{\beta_1}
\ar[ur]_{\alpha_1} &
}
\]

\begin{center}
\begin{figure}[H]
\[
\xymatrix@R-5pt
{
\K :
&
{s}
\ar@<1ex>[r]^{\alpha}
\ar@<-1ex>[r]_{\beta}
&
{t}
}
\]
\end{figure}
\end{center}
and the three following coverings $F_0$, $F_1$ and $F_2:\C\longrightarrow \K$ given by
$$F_i(s_0)=F_i(s_1)=s,\ \
F_i(t_0)=F_i(t_1)=t, \ \ F_i(\beta_0)=F_i(\beta_1)=\beta, \mbox{\ for\ } i=0,1,2$$
while
\begin{itemize}
\item $F_0(\alpha_0)=F_0(\alpha_1)=\alpha$
\item $F_1(\alpha_0)=F_1(\alpha_1)=\alpha + \beta$
\item $F_2(\alpha_0) = \alpha + \beta,\ F_2(\alpha_1) =\alpha$

\end{itemize}
Note that $F_0$ is the linearization of a covering of small categories, and $F_1$ is $F_0$ followed by an automorphism of $\K$. We will observe that $F_0$ and $F_1$ are in fact Galois coverings (see Definition \ref{defGalois}), while $F_2$ is not.
\end{exa}

\begin{defi}
A morphism from a covering
$F:\C\rightarrow\B$ to a covering $G:\D\rightarrow\B$ is a pair of $k$-linear functors $(H,J)$
where $H: \C  \to \D$, $J: \B \to \B$ are such that $J$ is an isomorphism, $J$ is the identity on objects and $GH=JF$.
 The category of coverings of $\B$
is denoted $\cov(\B)$.
\end{defi}

Our  next  purpose is to show that the automorphism group of a connected covering acts freely on
each fibre.

\begin{defi}\label{def.connected} A $k$-category $\B$ is connected if any two objects $b$
and $c$ of $\B$ can be linked by a finite walk made of non zero morphisms, more
precisely there exist a finite sequence of objects $x_1,\dots,x_n$ and non zero morphisms
$\varphi_1,\dots,\varphi_n$ such that $x_1=b,\ x_n=c$,  where $\varphi_i$ belongs either
to ${}_{x_{i+1}}\B_{x_i}$ or to ${}_{x_{i}}\B_{x_{i+1}}$.
\end{defi}

\begin{pro}\label{connected}
Let $F:\C\longrightarrow\B$ be a covering of $k$-categories. If $\C$ is connected, then
$\B$ is connected.
\end{pro}
\noindent\textbf{Proof. }\sf
Let $b$ and $c$ be objects in $\B_0$, and let $x_0$ and
$y_0$ be two objects respectively chosen in their fibres. Consider a walk of non zero
morphisms connecting $x_0$ and $y_0$ in $\C$. Since $F$ induces $k$-isomorphisms at each
star, the image by $F$ of a non zero morphism is a non zero morphism in $\B$. \qed

\begin{pro}\cite{le}\label{equal}
Let $F:\C\longrightarrow\B$ and $G:\D\longrightarrow\B$ be coverings of $k$-linear
categories. Assume $\C$ is connected. Two morphisms $(H_1,J)$, $(H_2,J)$ from $F$ to $G$ such that $H_1$ and $H_2$ coincide on some object are equal.
\end{pro}
\noindent\textbf{Proof. }\sf
Let $(H,J)$ be a morphism of coverings, let
$x_0$ be an object of $\C$ and consider the map between stars induced by $H$:
\[H^{x_0}_{H(x_0)}:\st_{x_0}\C\longrightarrow\st_{H(x_0)}\D.\]
Observe that $GH(x_0)=JF(x_0)$. There is a commutative diagram

\[
\xymatrix@R-5pt {\st_{x_0}\C
\ar[rr]^{H^{x_0}_{H(x_0)}}
\ar[d]_{F^{x_0}_{F(x_0)}}
&&
\st_{H(x_0)}\D
\ar[d]^{G^{H(x_0)}_{G(H(x_0))}}
\\
\st_{F(x_0)}\B \ar[rr]^{J^{F(x_0)}_{J(F(x_0))}} & & \st_{JF(x_0)}\B
}
\]

\noindent
where the morphisms $F^{x_0}_{F(x_0)}$, $J^{F(x_0)}_{J(F(x_0))}$ and $G^{H(x_0)}_{G(H(x_0))}$ are $k$-isomorphisms.
Consequently $H$ is an isomorphism at each star level, determined by $F$, $J$ and $G$. In case
$H_1$ and $H_2$ are morphisms such that $H_1(x_0)=H_2(x_0)$, the $k$-linear maps
$\st_{x_0}\C\longrightarrow\st_{H_1(x_0)}\D = \st_{H_2(x_0)}\D $ induced by $H_1$ and $H_2$ are
equal, hence $H_1$ and $H_2$ coincide on morphisms starting or ending at $x_0$, in particular
$H_1$ and $H_2$ coincide on objects related to $x_0$ by a non-zero morphism. Since $\C$ is
connected, it follows that $H_1$ and $H_2$ coincide on every object and on every morphism of
$\C$. \qed

\begin{cor}
Let $F: \C \to \B$ be a  connected covering of a $k$-linear category $\B$.  The group $\aut_1(F)=\{ (H,1): F \to F \mid H \ \mbox{an isomorphism of}\ \C\}$ acts freely on each fibre.
\end{cor}

\section{\sf Galois and universal coverings}\label{galoiscoverings}

We start this section with the definition of a Galois covering in order to study
properties of this kind of coverings.  The main results are the description of the
structure of Galois coverings, and the relation between Galois coverings and the group of
automorphisms. Finally we consider universal objects in the
category of Galois coverings of a fixed $k$-category $\B$.

\begin{defi}\label{defGalois}
A covering $F: \C\longrightarrow\B$ of $k$-categories is a \emph{Galois covering} if $\C$
is connected and if $\aut_1 F$ acts transitively on some fibre. We denote $\gal (\B)$ the full subcategory of $\cov (\B)$ whose objects are the Galois coverings of $\B$.
\end{defi}

It is natural to expect that $\aut_1 F$  should act transitively on each
fibre whenever it acts transitively on a particular one. In order to prove this fact, we
shall use a construction introduced in \cite{ga,boga}, see also \cite{cire}.

\begin{defi}
Let $G$ be a group acting by $k$-isomorphisms on a $k$-category $\C$, such that the
action on the objects is free, meaning that if $sx=x$ for some object, then $s=1$. The
set of objects of the categorical quotient $\C/G$ is the set of $G$-orbits of $\C_0$. The
$k$-module of morphisms from an orbit $\alpha$ to an orbit $\beta$ is
\[{ }_{_\beta}\left(\C/G\right)_\alpha= \left(\bigoplus_{x\in \alpha,\ y\in\beta}{}_y\C_x\right)/G\]
where for a $kG$-module $X$ we denote $X/G$ the $k$-module of coinvariants $X/(\ker\epsilon)X$, which is
the quotient of $X$ by the augmentation ideal, where $\epsilon :kG\rightarrow k$ is given
by $\epsilon(s)=1$ for all $s\in G$.
\end{defi}

\begin{rem}
The previous definition provides a $k$-category: the
composition is well defined precisely because   the action of $G$ on the objects is free.
\end{rem}

\begin{pro}
Let $G$ be a group acting by $k$-isomorphisms on a connected $k$-category $\C$, and
assume that the action on the objects is free. Then the projection functor
$P:\C\longrightarrow \C/G$ is a Galois covering with $\aut_1 (P)=G$.
\end{pro}

\noindent\textbf{Proof. }\sf
 The projection
functor is a covering since it is surjective on objects. For each choice of an object
$x_0\in \alpha$ and $y_0\in\beta$ we clearly have $k$-isomorphisms
\[\bigoplus_{y\in\beta}{}_y\C_{x_0}\rightarrow
{}_{_\beta}\left(\C/G\right)_\alpha \mbox{ \ and\  }
\bigoplus_{x\in\alpha}{}_{y_0}\C_{x}\rightarrow {}_{_\beta}\left(\C/G\right)_\alpha\]
which can be assembled in order to provide the required isomorphism of stars. Observe
that the fibres of $P$ are the orbits sets by construction, therefore the action of $\aut_1
P =G$ is transitive on each fibre.

Consider now $(H,1)\in\aut_1 P$ and let $x_0\in\C_0$. Since the action of $G$ on $\C_0$ is
free, there exists a unique $s\in G$ such that $sx_0=H(x_0)$. By definition, the element
$s$ provides an isomorphism of $\C$ such that $Ps=P$. The isomorphisms $(s,1), (H,1)$ of the
connected covering $P$ coincide on an object, consequently
 they are equal as isomorphisms of $P$ by Proposition \ref{equal}.\qed

\begin{lem}
Let $F:\C\longrightarrow\B$ be a connected covering of $k$-categories and suppose there
exists a singleton fibre. Then every fibre is a singleton and $F$ is an isomorphism of
$k$-categories.
\end{lem}
\noindent\textbf{Proof. }\sf
Let $b\in\B_0$ be an object such that
$F^{-1}(b)=\{x\}$. Since $\C$ is connected it is enough to show that for a non zero
morphism in $\C$ with target or source $x$, the other extreme object $y$ is such that
$F^{-1}(F(y))= \{ y\}$. We denote $c=F(y)$. Assume $\varphi\in {}_{y}\C_{x}$ is non-zero
and let $y'\in F^{-1}(c)$, then  $\varphi \in \st_x\C$ and $F(\varphi)\in\st_{b}\B$.
Moreover, $F(\varphi)$ belongs to $\st_c\B$. Since $F$ induces an isomorphism $F^{y'}_c:
\st_{y'}\C\longrightarrow\st_c\B$, there is a unique $k$-linear combination $\sum_z
{}_{y'}h_z$ of morphisms from the fibre of $b$ to $y'$ such that
$F^{y'}_c(\sum_z{}_{y'}h_z) = F(\varphi)$. Now the fibre of $b$ is reduced to $x$, which
means that there is a non-zero morphism $\psi\in{}_{y'}\C_x$ such that $F\psi=F\varphi$.
Note that $\psi$ also belongs to $\st_x\C$, and recall that $F_b^x$ is an isomorphism
between the corresponding stars. Hence $\varphi=\psi$ and in particular their ending
objects are the same, namely $y=y'$. Finally since all the fibres are singletons, the
star property of a covering implies immediately that $F$ is an isomorphism. \qed

We are now able to prove the following result.

\begin{pro}\label{projection} Let
$F:\C\longrightarrow\B$ be a Galois covering. Then $\aut_1 F$ acts transitively on each
fibre.
\end{pro}
\noindent\textbf{Proof. }\sf
First consider the categorical quotient $P:
\C\longrightarrow \C/\aut_1 F$. There is a unique functor $F': \C/\aut_1 F\longrightarrow \B$
such that $F'P=F$, defined as follows: let $\alpha$ be an object of $\C/\aut_1 F$, that is,
an orbit of $\C_0$ under the action of $\aut_1 F$. Choose an object $x\in\alpha$ and define
$F'\alpha=Fx$.  Clearly $F'$ is well defined on objects. In order to define $F'$ on
morphisms, let $\alpha$ and $\beta$ be objects in $\C/\aut_1 F$, and recall that
\[{ }_{_\beta}\left(\C/\aut_1 F\right)_\alpha=
\left(\bigoplus_{x\in \alpha,\ y\in\beta}{}_y\C_x\right)/\aut_1 F.\]
Next observe that the morphism
$$F: \bigoplus_{x\in \alpha,\ y\in\beta}{}_y\C_x  \to {}_{F'(\beta)}\B_{F'(\alpha)}$$
is sucht that $F(s \varphi)=F(\varphi)$ for any $\varphi \in \bigoplus_{x\in \alpha,\
y\in\beta}{}_y\C_x $ and any $s \in \aut_1 F$. Finally the commutative triangle of
morphisms between corresponding stars shows that $F'$ is indeed a
covering.\\
Since $F$ is a Galois covering, there exists a fibre where the action of $\aut_1 F$ is
transitive, which means that the corresponding fibre of $F'$ is a singleton. Since $F$ is
a Galois covering, $\C$ is connected as well as $\C/\aut_1 F$ by Proposition
\ref{connected}. The preceding Lemma asserts that all the fibres of $F'$ are singletons,
which exactly means that the action of $\aut_1 F$ is transitive on each fibre of $F$. \qed

As a consequence we obtain the following description of Galois coverings.

\begin{thm}\label{structure}
Let  $F:\C\longrightarrow\B$ be a Galois covering. Then there exists a unique
isomorphism of categories $F':\C/\aut_1 F\longrightarrow \B$ such that $F'P=F$, where
$P:\C\longrightarrow\C/\aut_1 F$ is the Galois covering given by the categorical quotient.
\end{thm}
\noindent\textbf{Proof. }\sf
The proof of the preceding results provides the
covering $F'$, which has a singleton fibre. Then all the fibres of $F'$ are singletons
and $F'$ is an isomorphism. \qed
\begin{exa} \cite{le,le2} An easy computation shows that $\aut_1(F_2)$ is trivial for the covering $F_2$ in Example \ref{nongalois}. However each fibre has two objects, hence the action of the trivial group is not
transitive on the fibres, consequently  $F_2$ is not Galois. Observe that $F_0$ and $F_1$ are Galois coverings.
\end{exa}
Next we recall a result of Patrick Le Meur concerning factorizations of Galois coverings.
\begin{lem}\label{Hsurjective}
Let $F:\C\longrightarrow \B$ and $G:\D\longrightarrow \B$ be Galois coverings, and let
$(H,J)$ be a morphism from $F$ to $G$. Then $H$ is surjective on objects.
\end{lem}
\noindent\textbf{Proof. }\sf
Let $H(c) = d\in\D_0$ be an object which is in the image of $H$. First we prove that any object $d'$ linked to $d$ by a non-zero morphism is also in the image of $H$. Let for instance $0\neq f\in {}_{d'}\D_d$. Note that $GH$ is a covering since $JF=GH$.
Considering $G(f)$, there exists a finite set of morphisms $(f_i)$ starting at $c$ and ending at objects $x_i$ such that
$GH(\sum f_i)=G(f)$, hence $G(\sum H(f_i))=G(f)$. Note that $H(f_i)$ is a morphism from $d$ to $H(x_i)$. Since $G$ is a covering
and $f$ is a morphism starting at $d$, we infer $f=\sum H(f_i)$. This implies that all the $H(x_i)$ coincide with $d'$,
hence $d'$ is in the image of $H$.

Finally, using that $\D$ is connected we conclude that any object of $\D$ is in the image of $H$.
\qed

\begin{thm} \cite{le,le2} \label{morphism2}
Let $F:\C\longrightarrow \B$ and $G:\D\longrightarrow \B$ be Galois coverings, and let
$(H,J)$ be a morphism from $F$ to $G$. Then there is a unique surjective group morphism $\Lambda :\aut_1
F\to\aut_1 G$ with $\Lambda(f,1)=(\lambda(f),1)$ such that $\lambda(f)H=Hf$ for each $(f,1)\in\aut_1 F$. Moreover $\ker\Lambda=\aut_1 H$ and $H$ is a Galois covering.
\end{thm}
\noindent\textbf{Proof. }\sf
Given $J$, we assert that the set of morphisms $(H',J)$ from $F$ to $G$
%verifying $GH'=FJ$
is in one-to-one correspondence with $\aut_1 G$ through the map which assigns $(gH,J)$ to each $g\in \aut_1 G$.
Firstly each $(gH,J)$ is a morphism from $F$ to $G$. Secondly if $(H',J)$ is such a morphism, given an object $c_0$ of $\C$,
both $H'(c_0)$ and $H(c_0)$ are in the same $G$-fibre. Since the action of $\aut_1 G$ is free and transitive on the fibres,
there exists a unique $g \in \aut_1 G$ such that $gH(c_0)=H'(c_0)$. Consequently $(H',J)$ and $(gH,J)$ are equal by Proposition \ref{equal}.

Then for each $(f,1)\in\aut_1 F$ there
exists a unique element $\Lambda(f,1) \in \aut_1 G$ such that $\lambda(f)H=Hf$. The uniqueness of $\Lambda(f)$ and the equalities $$\lambda(f_1f_2)H=Hf_1f_2=\lambda(f_1)Hf_2=\lambda(f_1)\lambda(f_2)H$$
imply that $\Lambda$ is a group morphism.

Moreover $\Lambda$ is surjective. Let
$(g,1) \in\aut_1 G$. Using the previous Lemma consider an object $c$ in the $H$-fibre of $gHc_0$. A simple computation shows that $c$ and $c_0$ are in
the same $F$-fibre. Since $F$ is a Galois covering, there exists $(f,1)\in\aut_1 F$ such that
$fc_0=c$, then $Hfc_0=Hc=gHc_0$. Then $Hf=gH$ and $\lambda(f)=g$.

Note that $(f,1)\in\ker\Lambda$ if and only if $Hf=H$ which means precisely that $(f,1)\in \aut_1
H$.

{In order to prove that $H$ is a Galois covering, we already know that $H$ is surjective on objects. The functor $H$ induces isomorphisms between stars since $GH=JF$, hence the same
equality is valid at the stars level where $F$, $J$ and $G$ induce $k$-isomorphisms. This proves that $H$ is a covering. In order to show that $H$ is Galois, let $x$ and
$x'$ be in the same $H$-fibre. They are also in the same $F$-fibre, hence there exists
$(f,1)\in\aut_1 F$ such that $fx=x'$. We assert that in fact $(f,1)\in\aut_1 H$: indeed, $(Hf,J)$ and $(H,J)$
are both morphisms from $F$ to $G$ with the same value on $x$, hence they are equal by
Proposition \ref{equal}.\qed

\begin{rem} Two isomorphic $k$-linear categories have isomorphic categories of Galois coverings.
\end{rem}

\begin{defi}
A \emph{universal covering} $U:\U \to \B$ is an object in $\gal(\B)$ such that for any Galois covering $F:\C \to \B$, and for any $u_0\in \U_0$, $c_0\in \C_0$ with $U(u_0)=F(c_0)$, there exists a unique morphism $(H,1)$ from $U$ to $F$ such that $H(u_0)=c_0$.    \end{defi}

In case of existence, a universal covering is unique up to isomorphisms of
Galois coverings. In general universal coverings do not exist, as the following Example shows.
It has been obtained by Geiss and de la Pe\~na in \cite{gepe}:

\begin{exa}
Let $k$ be a field and $char(k)=2$. Consider the $k$-linear categories

\[
\xymatrix
{ {\C}_1 : &  x_0 \ar[rr]^{\alpha_0} \ar[rrd]^{\beta_0}  & & y_0 \ar[rr]^{\gamma_0} \ar[rrd]^{\delta_0} & & z_0 \\
& x_1 \ar[rr]_{\alpha_1} \ar[rru]_{\beta_1} & & y_1 \ar[rr]_{\gamma_1} \ar[rru]_{\delta_1} & & z_1}
\]

\[
\xymatrix@R-5pt
{
\B : & x \ar@<1ex>[r]^{\alpha} \ar@<-1ex>[r]_{\beta} &
y \ar@<1ex>[r]^{\gamma} \ar@<-1ex>[r]_{\delta} & z
}
\]
with ${\C}_1$ satisfying all commutativity relations and ${\B}$
satisfying the relations
\[\gamma \alpha = \delta \beta , \qquad
\gamma \beta = \delta \alpha.
\]
It is clear that ${\C}_1$ is a Galois covering of ${\B}$. Since $\charac(k)=2$, if
we set $a=\alpha + \beta, b= \beta, c=\gamma + \delta, d= \delta$, we get that $\cal B$
satisfies the relations
\[ca=0, \qquad cb=da.\]
In this case,
\[
\xymatrix@R-5pt
{&&&& \vdots & \vdots \\
&&& x_{-1} \ar[r]^{b_{-1}}\ar[d]^{a_{-1}} & y_{-1} \ar[r]\ar[d]^{c_{-1}} & z_{-1}\\
{\C_2} : && x_{0} \ar[r]^{b_{0}}\ar[d]^{a_{0}} & y_{0} \ar[r]^{d_{0}}\ar[d]^{c_{0}} & z_{0} &\\
& x_{1} \ar[r]^{b_{1}} & y_{1} \ar[r]^{d_{1}} & z_{1} & &\\
&\vdots & \vdots &&&\\
}
\]
with all commutativity relations and $c_i a_{i-1}=0$, is also a Galois covering
of $\cal B$. Now ${\cal C}_1$ and ${\cal C}_2$ admit no proper Galois covering since they
are simply connected, see \cite{MP, AS},
and there is no morphism between them.\end{exa}

%Even in the case of coverings of quivers with relations, where the functors respect arrows and relations, there is no universal covering in general. An example has been obtained by Geiss and de la Pe\~na in \cite{gepe}.

Now we will study the Kronecker category $\K$.
Recall that this category is given by two objects $s, t$, one-dimensional morphism spaces $_{s}\K_s$ and $ _{t}\K_t$  while $\dim_k \ _{t}\K_s=2$ and $_{s}\K_t=0$. Observe that for each choice of a vector basis of $ _{t}\K_s$, the category $\K$ is presented by the quiver
\[
\xymatrix
{{s}
\ar@<1ex>[r]^{\alpha}
\ar@<-1ex>[r]_{\beta}
&
{t}
}
\]

We start with a description of all Galois coverings of the category $\K$ up to isomorphisms.

Let $\{a,b\}$ be a basis of $ _{t}\K_s$.  Let $\C_{\{a,b\}}$ be the free $k$-category presented by the quiver
\[ \xymatrix{
\vdots & \vdots\\
s_1 \ar[r]^{a_1} \ar[ur] & t_{1} \\
s_0 \ar[r]^{a_0} \ar[ur]^{b_0} & t_0 \\
\vdots \ar[ur] & \vdots
}\]
and let $F_{\{a,b\}}: \C_{\{a,b\}} \to \K$ be given by
\begin{itemize}
\item $F_{\{a,b\}}(s_i)=s$,
\item $F_{\{a,b\}}(t_i)=t$,
\item $F_{\{a,b\}}(a_i) = a$,
\item $F_{\{a,b\}}(b_i) = b$.
\end{itemize}

\begin{pro}
Let $\{a,b\}$ be a fix chosen basis of $ _{t}\K_s$.  Then any  Galois covering of $\K$ is isomorphic to $F_{\{a,b\}}$ if it has infinite Galois group, and to a quotient of it otherwise.
\end{pro}

\noindent\textbf{Proof. }\sf
It can be seen that any  Galois covering of $\K$ is isomorphic to $F_{\{c,d\}}$ or a quotient of it, where $\{c,d\}$ is a basis of $ _{t}\K_s$. Now  $F_{\{c,d\}}\simeq F_{\{a,b\}}$ with an isomorphism of type $(1,J)$.
\qed
Note that the Kronecker category has no universal covering since the Definition of a universal covering only takes into account morphisms of type $(H,1)$.

\section{\sf Fundamental group}\label{fundamental}

As quoted in the Introduction, our main purpose is to provide an intrinsic definition of
the fundamental group $\pi_1$ of a $k$-category, where $k$ is a commutative ring.
Previous definitions, provided for instance by J. A. de la Pe\~{n}a and R. Mart\'\i nez-Villa
\cite{MP}, see also K. Bongartz and P. Gabriel \cite{boga}, depend on the presentation of
the category as a quotient of a free $k$-category by an ideal generated by some set of
minimal relations. Different presentations of the same $k$-category may provide different
groups through this construction, see for instance \cite{asde,buca,le1}.

We will prove the following fact concerning the group $\pi_1$ that we will define: if the universal covering $U$ exists, then the group $\aut_1 U$ is isomorphic to
$\pi_1$; in this case any group obtained through the presentation construction is a
quotient of $\pi_1$.

\begin{defi}
Let $\B$ be a $k$-category, and let $b_0$ be a fixed object in $\B_0$. Consider $\gal(\B,b_0)$ the subcategory of $\gal\B$ with the same objects and morphisms $(H,J)$ with $J(b_0)=b_0$.
Let $\Phi: \gal (\B,b_0) \to \mathsf{Sets} $ be the \emph{fibre
functor} which associates to each Galois covering $F$ the $F$-fibre $F^{-1}(b_0)$. We
define
$$\pi_1(\B,b_0)=\Aut\Phi.$$
\end{defi}

%ojo que si cambiamos de definicion del universal, la accion de los automorfismos del fibra %en la fibra del universal no es transitiva!!

\begin{rem}\label{lindo}
This fundamental group $\pi_1(\B,b_0)$ is the group of natural isomorphisms $\sigma: \Phi
\to \Phi$. In other words an element of the fundamental group is a family of invertible
set maps $\sigma_F: F^{-1}(b_0) \to F^{-1}(b_0)$ for each Galois covering $F$, which are
compatible with morphisms of Galois coverings; namely for each morphism $(H,J): F\to G$ in
$\gal (\B, b_0)$ the corresponding square
\[
\xymatrix@R-1pt
{{F^{-1}(b_0)}
\ar[r]^{\sigma_F}
\ar[d]_{H}
&
{F^{-1}(b_0)}
\ar[d]^{H}
\\
{G^{-1}(b_0)}
\ar[r]^{\sigma_G}
&
{G^{-1}(b_0)}
}
\]
is commutative.
\end{rem}

In case the universal covering exists, our purpose is to prove that the fundamental group is isomorphic to its automorphism group.

\begin{pro}
Let $\B$ be a connected $k$-category and let $b_0$ be an object. Assume there exists a universal covering $U$. Then $\pi_1(\B,b_0)$ acts freely and transitively on $U^{-1}(b_0)$.
\end{pro}
\noindent\textbf{Proof.} \sf
Let $\sigma\in\pi_1(\B,b_0)$, we define the action by $$\sigma u = \sigma_U(u).$$
This is an action by the definition of composition of automorphisms of the fibre functor.

Let us first prove that the action is free. Assume $\sigma u =u$. Let $F$ be a Galois covering and let $c$ be an object in $F^{-1}(b_0)$. Consider  the unique morphism $(H,1)$ from $U$ to $F$ such that $Hu=c$. Using Remark \ref{lindo} we obtain $\sigma_F(c)=c$.

In order to prove transitivity, let $u$ and $u'$ be objects in $U^{-1}(b_0)$. We are going to define an automorphism  $\sigma$ of $\Phi$ such that $\sigma_U(u)=u'$. Let $F$ be a Galois covering and $c$ some element in $F^{-1}(b_0)$, and let $(H,1)$ be the unique morphism such that  $Hu=c$. We define $$\sigma_F(c)=H(u').$$ Using the uniqueness of the morphisms starting at the universal covering, one can prove that the family  $(\sigma_F)$ is indeed an automorphism of $\Phi$. \qed

\begin{pro}
The actions of $\pi_1(\B,b_0)$ and $\aut_1(U)$ on $U^{-1}(b_0)$ commute.
\end{pro}
\noindent\textbf{Proof.} \sf
This is an immediate consequence of Remark \ref{lindo}.\qed

We recall that an anti-morphism $\varphi:G\to G'$ is a map such that $\varphi(g_1g_2)=\varphi(g_2)\varphi(g_1)$
for any elements $g_1$ and $g_2$ in $G$. Of, course an anti-morphism $\varphi$ provides a unique usual group morphism $\psi$ given by $\psi(g)=\varphi(g^{-1})$.

\begin{lem}\label{commute}
Let $G$ and $G'$ be groups acting freely and transitively on a non empty set $X$. Assume the actions commute. Then each choice of an element in $X$ determines an anti-isomorphism from $G$ to $G'$.
\end{lem}
\noindent\textbf{Proof.} \sf
Choose an element $x\in X$. Define $\varphi : G\to G'$ by $gx =\varphi(g)x$. This map is well defined and bijective. Moreover it is a group {anti-morphism} precisely because the actions commute.  \qed

\begin{thm}
Suppose that a connected $k$-category $\B$ admits a universal covering $U$.
Then $$\pi_1(\B,b_0) \simeq \aut_1 U.$$
\end{thm}

\begin{cor}
Let $\B$ be a connected $k$-category admitting a universal covering $U$ and let $b_0$ and $b_1$ be two objects. Then
$\pi_1(\B,b_0)$ and $\pi_1(\B,b_1)$ are isomorphic.
\end{cor}

\begin{cor}
Let $\B$ be a $k$-category admitting a universal covering, and consider a presentation of
$\B$ given by a quiver $Q$ and an admissible two-sided ideal $I$ provided with a minimal
set of generators $R$ given by parallel paths. Let $\pi_1(Q,R,b_0)$ be the group of the
presentation as defined in \cite{MP}, with respect to a vertex $b_0$. Then there is a
group surjection $\pi_1(\B,b_0)\longrightarrow\pi_1(Q,R,b_0)$.
\end{cor}
\noindent\textbf{Proof. }\sf
In \cite{MP} it is proven that the group
$\pi_1(Q,R,b_0)$ can be realized as $\aut_1 F$ for a Galois covering $F$. Then the
universal covering $U$ of $\B$ provides an epimorphism from  $\aut_1 U$
to $\pi_1(Q,R,b_0)$.\qed

For any group $\Gamma$, let {$\Sigma_l(\Gamma)$} be the {full} subcategory {of the category of left $\Gamma$-sets} whose objects are sets with a transitive {left} action of the group $\Gamma$ such that the isotropy group of an element is invariant. Note that in this case the isotropy group of any element is invariant, since the action is transitive. {$\Sigma_r(\Gamma)$ denotes the analogous category, where objects are right $\Gamma$-sets.
}

We will prove that the fibre functor $\Phi$ is an equivalence when considered as a functor from $\gal_1(\B, b_0)$ to the category $\Sigma_l(\pi_1(\B,b_0))$, where $\gal_1(\B, b_0)$ is the subcategory of $\gal(\B, b_0)$ with same objects and morphisms of type $(H,1)$.

\begin{pro}
Let $\B$ be a $k$-category admitting a universal covering $U: \U \to \B$. Let {$$S: \gal_1(\B, b_0) \to \Sigma_r(\aut_1 U)$$} be the functor given by $S(F)=\Hom_{\gal_1(\B, b_0)}(U,F)$ and defined by composition on morphisms.  Then $S$ is an equivalence.
\end{pro}

\noindent\textbf{Proof. }\sf
First we assert that the action of $\aut_1(U)$ on $\Hom_{\gal_1(\B, b_0)}(U,F)$ is transitive.  Let $X,X' \in \Hom_{\gal_1(\B, b_0)}(U,F)$, let $t_0 \in \U_0$ and let $t_1 \in \U_0$ such that $X(t_0)=X'(t_1)$. Let $g$ be the unique endomorphism of $U$ such that $g(t_1)=t_0$, which exists since $U$ is universal. Moreover, $g$ belongs to $\aut_1 U$ with inverse given by the endomorphism $g'$ sending $t_1$ to $t_0$, and $Xg=X'$ since they coincide on $t_1$. Observe that in general the action is not free since $t_1$ is not uniquely determined.

Moreover this action has invariant isotropy group: let $X \in \Hom_{\gal_1(\B, b_0)}(U,F)$ and $g\in \aut_1 U$ such that $Xg=X$.  For any $h\in \aut_1 U$, using Theorem \ref{morphism2}, we have
$$ X h g h^{-1} = \lambda_X(h) X g h^{-1} = \lambda_X(h) X h^{-1} = X h h^{-1} = X.$$

Next we prove that $S$ is faithful. Let $F$ and $G$ be Galois coverings, and let $H, H' \in \Hom_{\gal_1(\B, b_0)}(F,G)$ such that $S(H)=S(H')$, that is, $HX=H'X$ for any $X \in \Hom_{\gal_1(\B, b_0)}(U,F)$, a non-empty set.  Then $H$ and $H'$ coincide on some object, and hence they are equal.

In order to prove that $S$ is full, let $f: \Hom_{\gal_1(\B, b_0)}(U,F) \to \Hom_{\gal_1(\B, b_0)}(U,G)$ be an $\aut_1 U$-morphism.  We are looking for a morphism $H \in \Hom_{\gal_1(\B, b_0)}(F,G)$ such that $S(H)=f$, that is, $f(X)=HX$ for any $X \in \Hom_{\gal_1(\B, b_0)}(U,F)$.  Given $X$ let $H$ be the unique morphism from $F$ to $G$ such that $H(X(t_0))= f(X)(t_0)$.  It is clear that $f(X)=HX$ since they coincide on some object.  For any other $X'$, we know that there exists $g \in \aut_1 U$ such that $Xg=X'$. Then
$$f(X')(t_0)=f(Xg)(t_0)=f(X)g(t_0)=f(X)(gt_0)$$
$$\lefteqn {=HX(gt_0)=HXg(t_0)=HX'(t_0)}$$
and hence $f(X')=HX'$.

Finally we prove that $S$ is dense.  Let $E \in \Sigma_r(\aut_1 U)$, $I$ the isotropy group of any element in $E$.  Let $P: \U \to \U/I$ be the projection functor, and let $F: \U/I \to \B$ as constructed in the proof of Proposition \ref{projection}. We define
$$\aut_1 U \ /I \to  \Hom_{\gal_1(\B, b_0)}(U,F) $$
given by $\overline g \mapsto Pg$.  This map is well-defined and injective.  In order to prove that it is surjective, let $H \in \Hom_{\gal_1(\B, b_0)}(U,F)$. Since $\aut_1 U$ acts transitively on $\Hom_{\gal_1(\B, b_0)}(U,F)$, there exists $g \in \aut_1 U$ such that $Pg=H$.
On the other hand, the map $\aut_1 U / I \to E$ given by $\overline{g} \mapsto ge$ is an isomorphism in
$\Sigma_r(\aut_1 U)$.
\qed

\begin{thm}
Let $\B$ be a $k$-category admitting a universal covering $U: \U \to \B$. Then the fibre functor $\Phi$ is an equivalence when considered as a functor from $\gal_1(\B, b_0)$ to the category {$\Sigma_l(\pi_1(\B,b_0))$}.
\end{thm}

\noindent\textbf{Proof. } \sf From the previous Proposition, it is enough to see that $\Phi$ and $S$ are naturally isomorphic.

Let $t_0$ be a fixed element in $U^{-1}(b_0)$ and let $\gamma:\pi_1(\B,b_0)\to \aut_1(U)$ be the {anti-isomorphism} defined by the equality $\sigma_U(t_0)= \gamma(\sigma)(t_0)$ (see Lemma \ref{commute}). The functor $\Sigma(\gamma)$ induced by $\gamma$ is a natural isomorphism from {$\Sigma_r(\Aut_1 U)$ to $\Sigma_l(\pi_1(\B,b_0))$.}

Given $F\in \gal_1(\B, b_0)$ one must prove that $\Sigma(\gamma)(\Hom_{\gal_1(\B, b_0)}(U,F))$ is naturally
isomorphic to $F^{-1}(b_0)$ as left $\pi_1$-sets. {Let $e$ be the map defined by $e(X)=X(t_0)$ for $X\in \Hom_{\gal_1(\B, b_0)}(U,F)$. Since $U$ is universal, $e$ is a bijection which is clearly natural. We assert that $e$ commutes with the left action of $\pi_1(\B, b_0)$. Let $\sigma\in\pi_1(\B,b_0)$, then $$e(\sigma X) = e(X \gamma(\sigma)) = \left(X\gamma(\sigma)\right)(t_0)=X(\gamma(\sigma)(t_0))=X\sigma_U(t_0).$$
On the other hand
$$\sigma e(X)= \sigma X(t_0) = \sigma_F (X(t_0).$$
Both elements are equal by means of Remark \ref{commute}.}
\qed

\section{\sf First Hochschild cohomology and Galois groups}

In this section our main purpose is to provide a canonical embedding from the additive
characters of the intrinsic $\pi_1$ that we have defined to the first Hochschild-Mitchell
cohomology vector space of $\B$. This will be achieved in case $k$ is a field and
assuming that the endomorphism  ring  of each object of the category is reduced to $k$,
{and that there exists a Galois covering whose group is isomorphic to the fundamental group (for instance if there exists a universal covering)}.

First we will provide an intrinsic and direct way of describing the injective morphism
from the additive characters of the group of automorphisms of a Galois covering to the
first Hochschild-Mitchell cohomology vector space of $\B$.

Assem and de la Pe\~na have described this map in \cite{asde} when $G$ is the fundamental
group of a triangular finite dimensional algebra presented by a quiver with relations. In
\cite{jap-s} de la Pe\~na and Saor\'\i n noticed that the triangular hypothesis is
superfluous. This map has also been obtained in a spectral sequence context in
\cite{cire}.

In order to provide the canonical morphism, we first recall (see \cite{gr,cima}) that a
Galois covering of $\B$ provides a  grading  for each choice of objects in the fibres.
As expected another choice of objects provides a conjugated grading. We will translate in
this setting the connectivity hypothesis of the Galois covering. Finally the definition
of Hochschild-Mitchell derivations as well as of  the inner ones will provide the context for
a natural definition of the required map.

\begin{defi}
Let $\B$ be a $k$-category, where $k$ is a ring. A $G$- grading  $Z$ of $\B$ by a group
$G$ is a decomposition of each  $k$-module of morphisms  as a direct sum of
 $k$-modules   $ Z_s$
indexed by $G$. For each couple of objects $b$ and $c$ we have ${}_c\B_b=\oplus_{s\in
G}Z_s\left({}_c\B_b\right)$ and
$$Z_t\left({}_d\B_c\right) Z_s\left({}_c\B_b\right) \subset Z_{ts}\left({}_d\B_b\right),$$
where elements of $Z_s\left({}_c\B_b\right)$  are called homogeneous morphisms of degree
$s$ from $b$ to $c$.
\end{defi}
The following result is clear:
\begin{pro}
Let $\B$ be a $k$-category as above, equipped with a $G$-grading  $Z$. Let
$\left(t_b\right)_{b\in\B_0}$ be a family of elements of $G$ associated to the objects of
$\B$. Define $Y_s\left({}_c\B_b\right)=Z_{t_cs{t_b}^{-1}}\left({}_c\B_b\right)$. Then $Y$
is also a $G$-grading  of $\B$.
\end{pro}

\begin{thm}\cite{cima}\label{graduation}
Let $\C\to\B$ be a Galois covering of categories. Let $\left(x_b\right)_{b\in\B_0}$ be a
choice of objects of $\C$, where $x_b$ belongs to the $F$-fibre of $b$ for each
$b\in\B_0$. Then there is an $\aut_1 F$-grading  of $\B$. Another choice of fibre objects
provides a grading  with the same homogeneous components as described in the preceding
Proposition.
\end{thm}
\noindent\textbf{Proof. }\sf
Since $F$ is a covering, each morphism space of $\B$ is equipped with a  direct sum
decomposition $\oplus_{y\in F^{-1}c}F\left({}_y\C_{x_b}\right)$. Moreover since $F$ is
Galois, for each $y$ in the $F$-fibre of $c$ there exists a unique automorphism $s$ such
that $y=sx_c$. This element $s$ will provide the degree of the direct summand, more
precisely
$$Z_s\left({}_c\B_b\right)=F\left({}_{sx_c}\C_{x_b}\right).$$
It is straightforward to check that this is indeed an $\aut_1 F$-grading. Moreover, a
different choice of objects in the fibres $\left(t_bx_b\right)_{b\in\B_0}$ where
$\left(t_b\right)_{b\in\B_0}$ is a family of elements of $\aut_1 F$ provides another
grading, which is precisely the grading
described in the proposition above. \qed

\begin{rem}
In \cite{cima} the converse is obtained: in case $\B$ is $G$-graded where $G$ is an
arbitrary group, the smash product  category construction provides a Galois covering with
automorphism group $G$. As expected, one recovers the original grading  as the
one  induced by the smash Galois covering as defined in \cite{cima}.
\end{rem}

\begin{defi}A \emph{homogeneous walk} $w$ in a $G$-graded $k$-category from an object $b$ to an
object $c$ is a sequence of non zero homogeneous paths. It consists of a sequence of
objects $x_1=b,\dots,x_i\dots,x_n=c$, a sequence of signs $\epsilon_1,\dots,\epsilon_n$
where $\epsilon_i\in\{-1,+1\}$, and non-zero homogeneous morphisms
$\varphi_1,\dots,\varphi_n$ such that if $\epsilon_i=1$ then $\varphi_i\in
{}_{x_{i+1}}\B_{x_i}$ while if $\epsilon_i=-1$ then $\varphi_i\in {}_{x_{i}}\B_{x_i+1}$.
The \emph{degree} of $w$ is the following ordered product of elements of $G$:
$$\deg w= (\deg\varphi_n)^{\epsilon_n}\cdots(\deg\varphi_i)^{\epsilon_i}\cdots(\deg\varphi_1)^{\epsilon_1}$$
\end{defi}

\begin{rem}
Note that a homogeneous non zero endomorphism involved in a homogeneous walk at
position $i$ appears  with its degree, or the inverse of its degree, according to the
value of $\epsilon_i$.
\end{rem}

\begin{defi} Let $\B$ be a $G$-graded $k$-category.
The grading is called \emph{connected} if for any couple of objects $b$ and $c$ of $\B$
and for any element $s\in G$, there exists a homogeneous walk from $b$ to $c$ of degree
$s$.
\end{defi}

\begin{thm}\label{connected.grading}
Let $F:\C\to\B$ be a Galois covering. Then the induced grading on $\B$ is connected.
\end{thm}
\noindent\textbf{Proof. }\sf
Let $\left(x_b\right)$ be a choice of an object in each $F$-fibre, providing a grading
of $\B$. Note that a morphism $\varphi$ in $\C$ from $sx_b$ to $tx_c$ has homogeneous
image of degree $s^{-1}t$, since $F\varphi=Fs^{-1}\varphi$. This observation shows that a
walk in $\C$ from $x_b$ to some $sx_d$ projects to a homogeneous walk from $b$ to $d$ of
degree $s$. Since $\C$ is connected, the theorem is proved.\qed

Let us briefly recall the definition of Hochschild-Mitchell cohomology in degree one (see
for instance \cite{mi}). This cohomology coincides with usual Hochschild cohomology of
algebras in case the $k$-category has a finite number of objects.

\begin{defi}
Let $\B$ be a $k$-category. A \emph{derivation} of $\B$ is a collection $D =
\left({}_cD_b\right)_{c,b\in\D_0}$ of $k$-linear endomorphisms of each $k$-module of
morphisms ${}_c\B_b$, such that $D(gf)=gD(f) + D(g)f$. More precisely if ${}_dg_c$ and
${}_cf_b$ are morphisms of $\B$, then $${}_dD_b(gf)=g{}_cD_b(f) + {}_dD_b(g)f.$$ An
\emph{inner derivation} $D_\alpha$ associated to a collection $\alpha$ of endomorphisms
at each object $\left(\alpha_b\right)_{a\in\B_0}$ is obtained in the usual way, namely
$$D_\alpha\left({}_cf_b\right)=\alpha_cf-f\alpha_b.$$
The \emph{first Hochschild cohomology} $k$-module $H^1(\B,\B)$ is the quotient of the
$k$-module of derivations by the inner ones.

\end{defi}

\begin{thm}
Let $k$ be a field and let $\B$ be a $k$-category such that the endomorphism ring of each
object is reduced to $k$. Let $F:\C\to\B$ be a Galois covering. Then there exists a
canonical injective morphism $$\Delta : \Hom(\aut_1 F, k^+)\longrightarrow H^1(\B,\B).$$
\end{thm}

\noindent\textbf{Proof. }\sf
Let $\chi:\aut_1 F\to k^+$ be an abelian character of $\aut_1 F$. In order to define
$\Delta\chi$ as a $k$-endomorphism of each $k$-module of morphisms of $\B$, we define
$\Delta_\chi$ on the homogeneous components of a grading  induced by the covering $F$.
Let $f$ be a morphism of degree $s$. By definition $\Delta\chi(f)=\chi(s)f$.

A standard computation shows that $\Delta\chi$ is indeed a derivation, which corresponds
to the well known construction of Euler derivations, see also \cite{fagegrma,fagrma}.
Moreover if the grading  is changed through a different choice of objects in the fibres
according to Theorem \ref{graduation}, the derivation $\Delta\chi$ is modified by an
inner derivation, hence the morphism $\Delta$ is canonic.

Assume $\Delta\chi$ is an inner derivation. Let $f$ be a non zero homogeneous morphism of
degree $s$, then $\chi(s)f=\alpha_cf-f\alpha_b=\left(\alpha_c-\alpha_b\right)f$ since
$\alpha_b\in k$ for each object $b$. Since $k$ is a field and $f$ is a non zero element
of a vector space,  $\chi(s)=\alpha_c-\alpha_b$. Moreover if $w$ is a homogeneous walk
of degree $s$ from $b$ to $c$, we also have $\chi(s)=\alpha_c-\alpha_b$ by an easy
computation. Finally since $F$ is Galois, we know by Theorem \ref{connected.grading} that
for each automorphism $s$ there exists a homogeneous walk from an object $b$ to itself of
degree $s$. Consequently $\chi(s)=0$ for every $s$ and $\Delta$ is injective.\qed

We have already observed that the Kronecker category $\K$ does not admit a universal covering, nevertheless there exist Galois coverings of $\K$ such that their group of automorphisms are isomorphic to the fundamental group. This observation provides interest to the following result.

\begin{cor}
Let $k$ be a field and let $\B$ be a $k$-category such that the endomorphism ring of each
object is reduced to $k$. Assume that $\B$ admits a Galois covering whose group is isomorphic to $\pi_1$ (for instance if $\B$ admits a universal covering).  Then there exists a
canonical injective morphism $$\Delta : \Hom(\pi_1(\B,b_0), k^+)\longrightarrow H^1(\B,\B).$$
\end{cor}

Observe that for the Kronecker category $\dim H^1(\K,\K)= 3$, then in general the above morphism is not an isomorphism.
%%%%%%%%%%%%%%%%%%%%%%%%%%%%%%%%%%%%%%%%%%%%%%%%%%%%%%%%%%%%%%%%%%%%%%%%%%%%%%%%%%%%%%%%%%%

%%%%%%%%%%%%%%%%%%%%%%%%%%%%%%%%%%%%%%%%%%%%%%%%%%%%%%%%%%%%%%%%%%%%%%%%%%%%%%%%%%%%%%%%%%%

\footnotesize \noindent C.C.:
\\Institut de Math\'{e}matiques et de Mod\'{e}lisation de Montpellier I3M, UMR 5149\\
Universit\'{e}  Montpellier 2, F-34095 Montpellier cedex 5,
France.\\
{\tt Claude.Cibils@math.univ-montp2.fr}

\noindent M.J.R.:
\\Departamento de Matem\'atica,
Universidad Nacional del Sur,\\Av. Alem 1253\\8000, Bah\'\i a Blanca, Argentina.\\ {\tt
mredondo@criba.edu.ar}

\noindent A.S.:
\\Departamento de Matem\'atica,
 Facultad de Ciencias Exactas y Naturales,\\
 Universidad de Buenos Aires,
\\Ciudad Universitaria, Pabell\'on 1\\
1428, Buenos Aires, Argentina. \\{\tt asolotar@dm.uba.ar}

\end{document}